\documentclass[11pt]{amsart}

\usepackage{amsmath}
\usepackage{fullpage}
\usepackage{xspace}
\usepackage[psamsfonts]{amssymb}
\usepackage[latin1]{inputenc}
\usepackage{graphicx,color}
\usepackage[curve]{xypic} 
\usepackage{hyperref}
\usepackage{graphicx}
\usepackage[greek,english]{babel}

\usepackage{amsmath}%
\usepackage{amsthm}%
\usepackage{amscd}
\usepackage{amsfonts}%
\usepackage{amssymb}%
\usepackage{graphicx}

\usepackage{mathrsfs}

\usepackage{tikz}
\usetikzlibrary{matrix,arrows}

\usepackage{tikz-cd}

\newtheorem{theorem}{Theorem}[section]

\newtheorem{conjecture}[theorem]{Conjecture}

\newtheorem{corollary}[theorem]{Corollary}

\newtheorem{lemma}[theorem]{Lemma}

\theoremstyle{remark}

\numberwithin{equation}{section}

\newcommand{\pfrak}{\mathfrak{p}}
\newcommand{\Pfrak}{\mathfrak{P}}
\newcommand{\qfrak}{\mathfrak{q}}

\newcommand{\Acal}{\mathscr{A}}

\newcommand{\Pro}{\mathbb{P}}

\newcommand{\Z}{\mathbb{Z}}

\newcommand{\Q}{\mathbb{Q}}

\newcommand{\rk}{\mathrm{rank}\,}

\newcommand{\alg}{\mathrm{alg}}

\makeatletter
\@namedef{subjclassname@2020}{%
  \textup{2020} Mathematics Subject Classification}
\makeatother

  \DeclareFontFamily{U}{wncy}{}
    \DeclareFontShape{U}{wncy}{m}{n}{<->wncyr10}{}
    \DeclareSymbolFont{mcy}{U}{wncy}{m}{n}
    \DeclareMathSymbol{\Sha}{\mathord}{mcy}{"58}

\begin{document}
\title[]{On Hilbert's tenth problem for number fields}

\author{Hector Pasten}
\address{ Departamento de Matem\'aticas,
Pontificia Universidad Cat\'olica de Chile, and CENIA. Chile}
\email[H. Pasten]{hector.pasten@uc.cl}%

\author{Xavier Vidaux}
\address{ Departamento de Matem\'aticas,
Universidad de Concepci\'on. Chile}
\email[X. Vidaux]{xvidaux@udec.cl}%


\thanks{H.P. was supported by ANID Fondecyt Regular grant 1230507 from Chile. X.V. was supported by Proyecto 2025001477INV from Universidad de Concepci\'on.}
\date{\today}
\subjclass[2020]{Primary: 11U05; Secondary: 11G05, 11J97} %
\keywords{Diophantine sets, number fields, Hilbert's tenth problem}%

\begin{abstract} Given a quadratic number field $K$ we show that the $abc$ conjecture over $K$ implies that $\mathbb{Z}$ is Diophantine in $K$ and, hence, a negative answer to Hilbert's tenth problem over $K$. More generally, the same conclusions hold for an arbitrary number field $K$ provided that the $abc$ conjecture holds over $K$ and that there is an elliptic curve over $\mathbb{Q}$ whose rank is positive and does not increase over $K$.
\end{abstract}

\maketitle



\section{Introduction}
Hilbert's tenth problem asked for an algorithm to decide solvability in $\Z$ of Diophantine equations over $\Z$. This was solved negatively by Davis--Putnam--Robinson and Matiyasevich \cite{DPR,Matiyasevich}. Of the many generalizations to other rings, the case of $\Q$ remained elusive for a long time. Recently, OpenAI proposed a negative solution in \cite{OpenAI1} (this depends on \cite{OpenAI2} and \cite{OpenAI3}) but, to the best of our knowledge, it is yet to be verified. Nevertheless, it led to a Diophantine definition of $\Z$ in $\Q$ by other researchers: a first AI-generated version by Glazer \cite{Glazer} and a considerably more accessible AI-assisted version by Burungale \cite{Burungale}. It is the authors's opinion that one can safely assume that \cite{Burungale} is complete; \cite{Glazer} is undergoing additional verifications by several researchers. This establishes the negative solution of Hilbert's tenth problem for $\Q$ in a strong form. See also \cite{KPV} where a very short argument is given, conditional on three standard arithmetic conjectures.

Our goal is to investigate Hilbert's tenth problem for number fields in general, not just $\Q$. We prove:

\begin{theorem}[General result]\label{ThmMain} Let $K$ be a number field and assume:
\begin{itemize}
\item[(i)] The $abc$ conjecture over $K$.
\item[(ii)] There is an elliptic curve $E$ over $\Q$ with $0<\rk E(\Q) = \rk E(K)$.
\end{itemize}
Then $\Z$ is Diophantine in $K$. In particular, Hilbert's tenth problem over $K$ has a negative solution.
\end{theorem}
From Zywina's work \cite{Zywina} we see that the required elliptic curves always exist for quadratic extensions, thus we deduce
\begin{corollary}[Quadratic case] Let $K$ be a quadratic number field and assume the $abc$ conjecture over $K$. Then $\Z$ is Diophantine in $K$. In particular, Hilbert's tenth problem over $K$ has a negative solution.
\end{corollary}

The $abc$ conjecture over $K$ is nothing but Vojta's conjecture on $\Pro^1$ with truncated counting functions \cite{Vojta} in the special case of the target divisor $[0]+[1]+[\infty]$. It implies the case of general target divisors by an argument of Granville using Belyi maps \cite{Granville} which originates in work of Elkies \cite{Elkies}. See Section \ref{SecHeights} for details.

We should note that there is the closely related question of Hilbert's tenth problem for \emph{rings of integers} of number fields, conditionally solved in \cite{MR, MP, Pasten} and, finally, unconditionally solved in \cite{KP, ABHS}. This is different; we are considering the whole number field, not its ring of integers. 

If we relax the condition that the definition of $\Z$ in $K$ be Diophantine, we can mention \cite{Robinson, PoonenDef, Koenigsmann, Park}. Note that Park's \cite{Park} concerns $O_K$ in $K$, but now we know that $\Z$ is Diophantine in $O_K$ by \cite{KP, ABHS}. Our route is different: under suitable assumptions we show that $\Q$ is Diophantine in $K$ and then we use the recent result that $\Z$ is Diophantine in $\Q$. Note that $\Q$ Diophantine in $K$ for $K\ne \Q$ would show that the arithmetic analogue of Corollary 2 in Koll\'ar's paper \cite{Kollar} is false.


\section{Heights and counting functions} \label{SecHeights}

Let $K$ be a number field. We write $h(-)$ for the standard absolute logarithmic height on $\Pro^1$, normalized by degree. For a divisor $D$ on $\Pro^1$ defined over $K$ we let $N^{(1)}(D,-)$ be the truncated counting function, in the sense of \cite{Vojta}, also normalized by degree. If $G\in K[T]$ is squarefree it defines a reduced effective divisor $D_G$ over $K$ on the projective line; for $x\in K$ not a zero of $G$ we write
$$
N^{(1)}(G,x) = N^{(1)}(D_G,x).
$$
Explicitly:
$$
N^{(1)}(G,x) = \frac{1}{[K:\Q]}\sum_{v_\pfrak(G(x))>0} \log {\rm N}\pfrak + O(1).
$$
For $x\in K^\times$ we define 
$$
N'(x) = \frac{1}{[K:\Q]}\sum_{v_\pfrak(x)>0 \,\, {\rm odd}} \log {\rm N}\pfrak.
$$
This is like the truncated counting function, but only considering primes with odd positive valuation.

Let us recall Vojta's conjecture over $K$ (see \cite{Vojta}):
\begin{conjecture}[$V_K$] Let $D$ be an effective reduced divisor on $\Pro^1$ defined over $K$, of degree $r$. Let $\epsilon>0$. Then we have
$$
(r-2-\epsilon)h(x) < N^{(1)}(D,x) + O_\epsilon(1)
$$
as $x$ varies in $\Pro^1(K)$ not in the support of $D$.
\end{conjecture}

The conjecture $V_\Q$ is equivalent to the classical $abc$ conjecture via a construction using Belyi maps, see \cite{Granville}. Over a number field $K$, the case $D=[0]+[1]+[\infty]$ is the direct analogue of the $abc$ conjecture and, in the same way, by using Belyi maps one sees that it implies the general case of $V_K$. Thus, we will work under the conjecture $V_K$.

From now on we fix \emph{one} squarefree polynomial $G\in \Z[T]$ (not just in $K[T]$) of degree $d\ge 7$. We only need to assume $V_K$ for the divisor $D_G$ defined by $G$.

\begin{lemma}[Lower bound for odd valuations]\label{LemmaVojta} Assume $V_K$ and let $\epsilon>0$. Then for all $x\in K$ with $G(x)\ne 0$ we have
$$
(d-4-\epsilon)h(x)\le N'(G(x)) + O_\epsilon(1).
$$
\end{lemma}
\begin{proof} Let $\gamma(n)$ be $0$ or $1$ depending on whether $n\ge 1$ is even or odd. From the inequality $2\le n+\gamma(n)$ applied to $n=v_\pfrak(G(x))$ if positive, we deduce
$$
2N^{(1)}(G,x) \le h(G(x)) + N'(G(x)) + O(1) \le dh(x)+ N'(G(x)) + O(1).
$$
The result now follows from $V_K$.
\end{proof}

\section{Pheidas's predicates}

A \emph{Pheidas predicate} for $K$ is a unary predicate $\Phi\subseteq K$ which is Diophantine over $K$ and satisfies the following conditions:
\begin{itemize}
\item[(Ph1)] All rational numbers whose denominator in reduced form is a square belong to $\Phi$.
\item[(Ph2)] There is a finite set of primes $S_\Phi$ of $K$ such that for every $\pfrak\notin S_\Phi$ and for every $b\in \Phi$ we have that $v_\pfrak(b)$ is not a negative odd integer.
\end{itemize}

Even-order poles predicates in $\Q$ were first introduced by Pheidas in \cite{Pheidas} as part of an approach to Hilbert's tenth problem for $\Q$. See also \cite{CornelissenZahidi, DemeyervanGeel}.

Using recent work of Burungale \cite{Burungale} we have

\begin{theorem} Every number field admits a Pheidas predicate.
\end{theorem}
\begin{proof} Let $\phi(b)$ be the positive existential unary formula from Definition 3.1 in \cite{Burungale} and take $\Phi\subseteq K$ as the set defined by it over $K$.

By Proposition 3.14 in \cite{Burungale} and using the same witnesses from $\Q$, we get condition (Ph1). Condition (Ph2) is a direct elementary verification with valuations, identical to the one in Proposition 3.2 of \cite{Burungale}.
\end{proof}
From now on we fix a Pheidas predicate $\Phi$ for $K$.


\section{Elliptic curves}

From now on we assume that there is an elliptic curve $E$ over $\Q$ with
$$
0<\rk E(\Q) = \rk E(K).
$$
Let us fix such an elliptic curve $E$ and let $m$ be a positive integer with $mE(K)\subseteq E(\Q)$. We write $\Gamma = mE(K)$. 

Let us fix a minimal Weierstrass equation for $E$ over $\Z$ and let $x,y$ be the coordinate maps on $E$. Let $\xi=-x/y$; this is a local parameter for $E$ at its identity element $O$. Define
$$
\Gamma^\circ = \{R\in \Gamma -\{O\}: x(R)y(R)\ne 0\}
$$
and
$$
\Acal = \{\xi(U)/\xi(V) : U,V\in \Gamma^\circ\}.
$$
Note that these sets are Diophantine over $K$, but nevertheless $\Gamma^\circ \subseteq E(\Q)$ and $\Acal \subseteq \Q^\times$.

Up to minor modifications, the following result appears in \cite{PoonenE} and we include a proof for the convenience of the reader.
\begin{lemma}[Elliptic curve approximations]\label{LemmaE} Let $k$ and $M$ be integers with $M\ge 2$. There is $\beta\in \Acal$ depending on them, such that $\beta-k \in M\Z_p$ for each $p|M$.
\end{lemma}
\begin{proof}
We may assume $k\ne 0$ (otherwise use $k=M$). 

Let us fix a non-torsion point $P\in\Gamma$. For each prime $p | M$, the formal group of the minimal Weierstrass equation that we fixed provides an open subgroup $H_p\subseteq E(\Q_p)$ of finite index with local parameter $\xi$, such that $\xi(H_p)\subseteq M\Z_p$ and $xy\ne 0$ on $H_p\setminus\{O\}$. 

Let $s\ge 1$ be divisible by all the indices $[E(\Q_p):H_p]$ for $p|M$, and set $V=[s]P$ and $U=[k]V$. They are non-torsion points of
$\Gamma$ that belong to each $H_p$ and in particular $U,V\in\Gamma^\circ$. In the formal group at $p$ the multiplication-by-$k$ map has the form $[k](z)=kz+z^2F_k(z)$ with $F_k\in\Z[[z]]$, and we finally deduce that for $p|M$ we have
$$
 \frac{\xi(U)}{\xi(V)} - k = \xi(V) F_k\left(\xi(V)\right)\in M\Z_p.
$$
\end{proof}


\section{A predicate selecting enough rational numbers}

Let us introduce the Diophantine predicate
$$
\Theta(t):\qquad G(t)\ne 0 \wedge \left(\exists \lambda,\beta\in \Acal,\quad \frac{(\lambda-G(t))^2}{G(t)^3}\in \Phi\quad\wedge\quad \frac{(\beta-t)^2}{\lambda}\in \Phi\right).
$$

The two membership conditions in this definition will be called \emph{first and second memberships} respectively. Our goal is to eventually show (assuming the existence of $E$ and that $V_K$ holds) that $\Theta$ defines in $K$ a subset of $\Q$ containing $\Z$, up to a finite error. From there it will follow that $\Q$ (hence, $\Z$) is Diophantine in $K$.

For the moment we have:

\begin{lemma}[Integers satisfy $\Theta$]\label{LemmaZinTheta} For every $n\in\Z$ with $G(n)\ne 0$ we have $K\models \Theta(n)$.
\end{lemma}
\begin{proof}
Write $g=G(n)$. Lemma \ref{LemmaE} with $(k,M)=(g,2g^2)$ gives a rational number $\lambda\in\Acal$ with $\lambda-g\in g^2\Z_p$ for each $p| g$. Thus, $(\lambda-g)^2/g^3$ is $p$-integral for each $p| g$ and has even valuation at all other primes; we conclude it belongs to $\Phi$.

Write $\lambda=A/B$ in lowest terms. Lemma \ref{LemmaE} with $(k,M)=(n,2|AB|)$ gives a rational number $\beta\in\Acal$ with $\beta-n\in AB\Z_p$ for each $p| AB$. We deduce that $(\beta-n)^2/\lambda$ is $p$-integral for each $p| AB$ and has even valuation at all other primes, so it belongs to $\Phi$.
\end{proof}

\section{Odd valuations}

From now on we let $\Sigma$ be the set of rational primes below primes in $S_\Phi$ or that ramify in $K$. Let $S$ be the set of primes of $K$ above the primes in $\Sigma$.

\begin{lemma}[Propagation of odd valuations]\label{LemmaSameV} Let $b\in K$ be such that $K\models \Theta(b)$. Let $\lambda,\beta$ be witnesses for $\Theta(b)$, and let $\pfrak\notin S$ be such that $v_\pfrak(G(b))$ is positive and odd. Then, for every prime $\qfrak$ of $K$ above the same rational prime as $\pfrak$ we have
$$
v_\qfrak(\lambda)=v_\pfrak(G(b))\qquad\mbox{and}\qquad v_\qfrak(\beta-b)>0 .
$$
\end{lemma}
\begin{proof}
Define $g=G(b)$ and $e=v_\pfrak(G(b))$. The first membership is $g^{-1}(\lambda/g-1)^2\in\Phi$. If $\lambda\ne g$, its valuation at $\pfrak$ is $2v_\pfrak(\lambda/g-1)-e$, which is odd, hence positive by (Ph2). So $v_\pfrak(\lambda/g-1)>0$ (possibly $+\infty$) in all cases, and $v_\pfrak(\lambda)=e$.

Let $p$ be the rational prime below $\pfrak$. As $p\notin\Sigma$, it is
unramified in $K$. Since $\lambda\in\Q$, we get $v_\qfrak(\lambda)=e$ for every $\qfrak | p$. At any given $\qfrak | p$, the expression in the second membership has valuation $2v_\qfrak(\beta-b)-e$ if $\beta\ne b$. This number is odd, hence positive by (Ph2), from which we conclude $v_\qfrak(\beta-b)>0$.
\end{proof}
The next lemma is our key to discarding irrational elements of $K$ via $\Theta$.
\begin{lemma}[Upper bound for odd valuations]\label{LemmaAut} If $K\models \Theta(b)$ and $\sigma:K\to\Q^\alg$ is an embedding with $\sigma(b)\ne b$, then
$$
N'(G(b))\le h(b-\sigma(b)) + O(1).
$$
\end{lemma}
\begin{proof}
Let $\lambda,\beta\in\Acal\subseteq \Q^\times$ be witnesses for $\Theta(b)$. Let $P(b)$ be the product of the rational primes $p\notin \Sigma$ for which there is a prime $\pfrak | p$ in $K$ with $v_\pfrak(G(b))$ positive and odd. We observe that
\begin{equation}\label{EqnBd1}
N'(G(b)) \le \log P(b) + O(1)
\end{equation}
where the bounded error term comes from the contribution of primes in $S$, which is bounded uniformly in $b$ because $N'(-)$ does not consider multiplicities.

By Lemma \ref{LemmaSameV} we have $v_\qfrak(\beta-b)>0$ at every $\qfrak | P(b)$ of $K$.  Let $L$ be the Galois closure of $K$ and extend $\sigma$ to an automorphism of $L$. We note that every $\Pfrak|P(b)$ is unramified and for such a prime there is $\qfrak|P(b)$ with $\Pfrak|\qfrak$, which gives 
$$
v_\Pfrak(\beta-b) = v_\qfrak(\beta-b)>0.
$$
As $\sigma$ permutes the primes $\Pfrak$ in $L$ above a rational prime $p|P(b)$ while it fixes $\beta\in \Q$, we obtain
$$
v_\Pfrak(\beta-\sigma(b))>0
$$
for each $\Pfrak|P(b)$. 

Recall that the element $\delta= b-\sigma(b)$ is non-zero by assumption. Since $\delta= (\beta-\sigma(b)) - (\beta-b)$ we obtain $v_\Pfrak(\delta)>0$ for each $\Pfrak|P(b)$. It follows that
$$
h(\delta)\ge \frac{1}{[L:\Q]} \sum_{\Pfrak | P(b)} v_\Pfrak(\delta) \log {\rm N}\Pfrak \ge \frac{1}{[L:\Q]} \sum_{\Pfrak | P(b)}  \log {\rm N}\Pfrak = \log P(b).
$$
We conclude by \eqref{EqnBd1}.
\end{proof}


\section{Proof of the main result}

\begin{proof}[Proof of Theorem \ref{ThmMain}] As before, we are working under the assumption of $V_K$ and the existence of the elliptic curve $E$.

Let $X\subseteq K$ be the set defined by $\Theta$ over $K$. By Lemma \ref{LemmaZinTheta} we know that $X$ contains all of $\Z$ except for a finite set $F\subseteq \Z$.

Now let $b\in X-\Q$. As $b$ is irrational, there is an embedding $\sigma: K\to \Q^\alg$ that moves $b$, namely,  $b-\sigma(b)\ne 0$. Since moreover $b\in X$, Lemma \ref{LemmaAut} gives
$$
N'(G(b))\le h(b-\sigma(b)) +O(1) \le h(b)+ h(\sigma(b)) + O(1) \le 2h(b) +O(1).
$$
Now, Lemma \ref{LemmaVojta} gives for every $\epsilon>0$
$$
(d-4-\epsilon)h(b)\le N'(G(b)) + O_\epsilon(1).
$$
As $d=\deg G\ge 7$ we obtain
$$
(3-\epsilon)h(b) \le  2h(b) +O_\epsilon(1).
$$
With $\epsilon=1/2$ (say) we deduce that $h(b)$ is bounded uniformly in $b$. Thus, by Northcott we get that $F'=X-\Q$ is finite.

The set $Y=(X\cup F) - F'$ is Diophantine in $K$, it contains $\Z$, and it is contained in $\Q$. Taking fractions we get that $\Q$ is Diophantine in $K$. As $\Z$ is Diophantine in $\Q$ by \cite{Glazer, Burungale}, we get that $\Z$ is Diophantine in $K$.
\end{proof}


\section{Acknowledgements and use of AI}

H.P. was supported by ANID Fondecyt Regular grant 1230507 from Chile. X.V. was supported by Proyecto 2025001477INV from Universidad de Concepci\'on.

Our starting point was our previous work \cite{KPV} with A. Kozachinskiy. Using those ideas we brainstormed with AI models, leading to a preliminary result for $\mathbb{Q}(i)$. Further AI-assisted work, with the main ideas coming from the authors, led to the current results.


\end{document}